# The Conics Generated by the Method of Application of Areas

*A Conceptual Reconstruction*

**Dimitris Sardelis and Theodoros Valahas**

**Abstract**

The method of application of areas as presented in Euclid's Elements, is employed to generate the three conics as the loci of points with cartesian coordinates satisfying the equations

$$x^2 = L\,y \pm \lambda\, y^2,$$

where $L$ and $\lambda$ are constants defined as the initial settings of the geometric constructions produced by the applications method. This conceptual reconstruction supports the view that the conics were most probably discovered as plane curves by fusing the method of the application of areas with the concept of locus, long before Apollonius studied them as conic sections.

## A Glimpse on a Fragmented History

There is a widespread belief that the three curves: the parabola, the ellipse and the hyperbola, were first perceived as sections of cones, henceforth their commonly used names as conic sections, or conics in short. Indeed, Menaechmus (see [1, Volume 2, p.p 110-116]) is said to have constructed the three conics by taking three types of circular cones, right-angled, acute-angled, and obtuse-angled, and by slicing each with a plane that does not pass through its vertex. It is also a fact that Apollonius [2] constructed the three conics in a more unified way as sections of any single circular cone, right or oblique. It must be noted, however, that these constructions provide no clue as to why and how the terms parabola, ellipse and hyperbola were introduced. Subsequently, one cannot rule out the possibility that these three curves -as well as the terms designated to them- were discovered prior to their construction as conics.

Commenting on a proposition in Euclid's Elements [4, II.44], Proclus [3, p.419] suggests that the parabola, the ellipse and the hyperbola were first conceived as plane curves generated by the Pythagorean method of application of areas and that their names were so assigned as to literally express the means of their construction. The passage runs as follows:

"These things, says Eudemus are ancient and are discoveries of the Muse of the Pythagoreans, I mean the application of areas ($\pi\alpha\rho\alpha\beta o\lambda\acute{\eta}$ $\chi\omega\rho\acute{\iota}\omega\nu$), their exceeding and their falling-short. It was from the Pythagoreans that later geometers took the names, which they again transferred to the so-called conic lines, designating one of these a parabola ($\pi\alpha\rho\alpha\beta o\lambda\acute{\eta}$), another a hyperbola ($\upsilon\pi\varepsilon\rho\beta o\lambda\acute{\eta}$) and another an ellipse ($\acute{\varepsilon}\lambda\lambda\varepsilon\iota\psi\iota\varsigma$), whereas those godlike men of old saw the things signified by these names in the construction, in a plane, of areas upon a finite straight line. For, when you have a straight line set out and lay the given area exactly alongside the whole of the straight line, then they say that you apply ($\pi\alpha\rho\alpha\beta\acute{\alpha}\lambda\lambda\varepsilon\iota\nu$) the said area; when however you make the length of the area greater than the straight line itself, it is said to exceed ($\upsilon\pi\varepsilon\rho\beta\acute{\alpha}\lambda\lambda\varepsilon\iota\nu$), and when you make it less, in which case, after the area has been drawn, there is some part of the straight line extending beyond it, it is said to fall short ($\varepsilon\lambda\lambda\varepsilon\acute{\iota}\pi\varepsilon\iota\nu$). Euclid too, in the sixth book, speaks in this way both of exceeding and falling-short; but in this place he needed the application simply, as he sought to apply to a given straight line an area equal to a given triangle in order that we might have in our power, not only the construction of a parallelogram equal to a given triangle, but also the application of it to a finite straight line...Such then is the application handed down from early times by the Pythagoreans."

Following Proclus's suggestion, the aim of the present article is to propose a conceptual reconstruction of the way the three conics were generated by the method of application of areas. Thus, based on a set of propositions stated and proved in Euclid's Elements, we formulate and solve three distinct geometrical-construction problems the key element of which is to construct a rectangle and a square of equal area both lying alongside a finite straight line segment with a common corner at one of its edges. These two simple figures intersect at a point opposite to the line segment. If the area of the rectangle is considered to vary continuously, so is that of the square and their intersection point traces out a conic. In particular, (1) when the rectangle is taken to apply exactly on the given line segment, i.e., the given line segment forms the rectangle's base, the curve produced is a parabola, (2) when the rectangle is taken to apply on part of the given line segment, i.e., the rectangle's base is shorter than the given line segment, the curve produced is an ellipse, and (3) when the rectangle is taken to apply extending beyond the given line segment, i.e., the rectangle's base is longer than the given line segment, the curve produced is a hyperbola.

In absence of direct historical documentation, one cannot be certain whether the proposed reconstruction-scenario resembles at all the actual story. However, indirect evidence suggests that conics in general and conics represented as loci ($\tau\acute{o}\pi o\iota$), in particular, had been thoroughly studied long before Apollonius. Thus, we are told by Pappus [5, Book VII,p. 636] of four books of Euclid's Conics and of five books of Solid Loci by Aristaeus, both works being extant even in Pappus's times. A second indirect evidence comes from Archimedes [6] who begins his treatise on the Quadrature of the Parabola by stating without proof "well known" propositions taken from what he calls "Elements of Conics", most probably referring to these two treatises. Finally, a third indirect evidence comes from Apollonius himself when in the introduction of his work (in eight books) he states explicitly that the first four books deal mostly with subjects already treated by others before him.

Summing up, the scenario proposed in this article of how the conics were discovered may appear too simple not to be true and also quite feasible not to be realistic. For, all that was needed was to fuse the method of application of areas ($\pi\alpha\rho\alpha\beta o\lambda\eta$ $\chi\omega\rho\iota\omega\nu$) with the concept of locus (τόπος).

## The Application of Areas in Euclid's Elements

The objective of the method of application of areas as presented in Euclid's Elements is to construct ("apply") a rectangle (or, more generally, a parallelogram) on a unit base-line (a given line segment) whose area equals that of any given rectilinear figure (exact application of areas). Thus, the given rectilinear figure is a triangle in proposition [**4**, I.44] and two adjacent triangles (i.e., a quadrilateral) in proposition [**4**, I.45]. The application effected by the latter proposition is the construction of two adjacent rectangles (or, more generally, two adjacent parallelograms) which, taken jointly, they form a rectangle (or, more generally, a parallelogram) of equal area to that of the given quadrilateral. Since it must have been obvious that every rectilinear figure can be resolved ultimately into triangles, the exposition of proposition [**4**, I.45] seems to establish that the application of areas is an additive operation and, consequently, that the method is applicable for any rectilinear figure of any area. The scheme of the exact application of areas is complemented with proposition [**4**, II.14] where a square is constructed on the given line that is similarly situated to the rectangle (or, more generally, a parallelogram) produced in proposition [**4**, I.45] and has the same area. This proposition might be rightly entitled as the quadrature (squaring) of any rectilinear figure.

The above scheme is then extended by two propositions to cases where the application is not exact : On a given line segment to apply a rectangle (or, more generally, a parallelogram) whose area equals that of a given rectilinear figure, (i) *falling short* [**4**, VI.28] or (ii) *exceeding* [**4**, VI.29] by a rectangle (or, more generally, a parallelogram) similar to a given rectangle (or, more generally, a parallelogram). These constructions could be called as *deficient* and *excessive* applications, respectively. No proposition analogue to [**4**, II.14] is found in Euclid's Elements for these cases. However, once the required rectangle has been constructed in each case, it is straightforward to also construct a similarly situated square of equal area.

The rectangle-square constructions effected by the (exact and inexact) application of areas method, are incorporated into Problems 1, 2 and 3 presented below. As for the key proposition [**4**, II.14] used for the subsequent derivation of the three conics (Corollaries 1, 2 and 3), we adopt Euclid's style of exposition and generalize its format for any given rectangle (with any length and width).

## The Parabola

   *Problem 1.* Suppose that to a given straight line segment there has been applied a rectangle equal to a given rectilinear figure. It is required to construct a square equal and similarly situated to the applied rectangle with a corner at one of the edges of the given straight line segment .

Let **AB** be the given straight line segment, **X** the given rectilinear figure and **ABCD** the rectangle applied on **AB** that equals **X** [**4**, I.45]. Thus, it is required to construct a square equal and similarly situated to **ABCD** with a corner at **A**.

If **AB** equals **BC**, then that which was proposed is done; for a square **ABCD** has been constructed equal to the rectilinear figure **X**.

But, if not, one of the straight line segments **AB** or **BC** is greater.

Let **AB** be greater than **BC** (see Figure 1), and extend **AB** to **E.** Make **AD** equal to **EA**, and bisect **EB** at **F** (by [**4**, I.3 and I.10]).

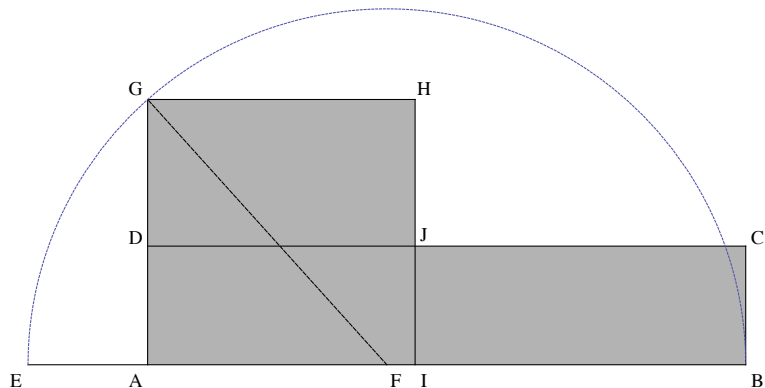

Figure 1

Describe the semicircle **EGB** with center **F** and radius one of the straight line segments **EF** or **FB.** Extend **AD** to **G**, and join **G** with **F** (by [**4**, I. Def. 18]).

Then, since the straight line segment **EB** has been cut into equal segments at **F** and into unequal segments at **A**, the rectangle **EA** by **AB** together with the square on **AF** equals the square on **EF** [**4**, II.5]:

$$(EA)\,(AB) + (AF)^2 = (EF)^2$$

But **EF** equals **FG**. Therefore the rectangle **EA** by **AB** together with the square on **AF** equals the square on **FG**:

$$(EA)\,(AB) + (AF)^2 = (FG)^2$$

But the sums of squares on **AF** and **AG** equals the square on **FG** [**4**, I.47] (the Pythagorean Theorem). Therefore, the rectangle **EA** by **AB** together with the square on **AF** equals the sum of the squares on **AF** and **EA**:

$$(EA)\,(AB) + (AF)^2 = (AF)^2 + (AG)^2$$

Subtract the square on **AF** from each. Therefore, the remaining rectangle **EA** by **AB** equals the square on **AG**, i.e, the square **AIHG**:

$$(EA)(AB) = (AG)^2$$

But the rectangle **EA** by **AB** is **ABCD**, for **EA** equals **AD.** Therefore the rectangle **ABCD** equals the square **AIHG**:

$$(AB)(BC) = (AG)^2 \qquad (1)$$

Therefore, if **AB** is greater than **BC**, a square **AIHG** has been constructed equal and similarly situated to the rectangle **ABCD**.

If **AB** is less than **BC**, then the same rationale also yields the construction of a square **AIHG** equal and similarly situated to the rectangle **ABCD** as displayed in Figure 2.

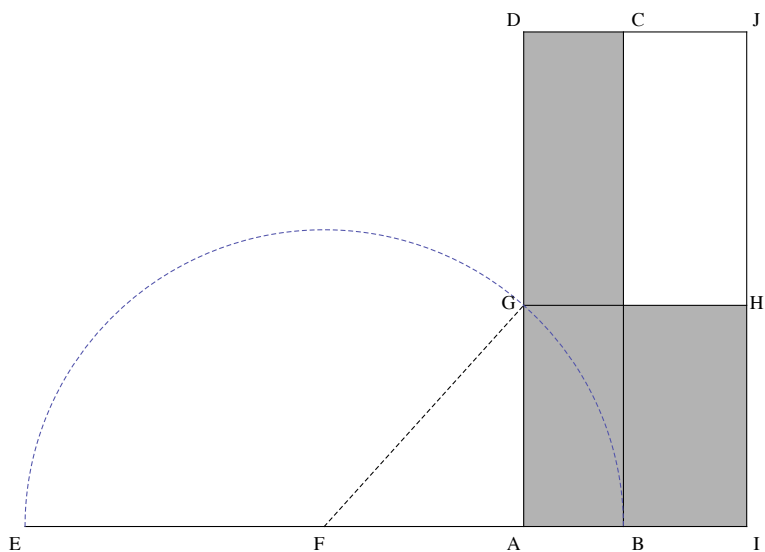

Figure 2

So, if **AB** equals **BC**, the constructed square applies to the entire line segment **AB**. If **AB** is greater than **BC**, the constructed square applies to part of the line segment **AB**. In either case, the boundaries of **AIHG** and of **ABCD** intersect at a point **J** (Figure 1). On the other hand, if **AB** is less than **BC**, the constructed square applies beyond **AB**, on the line segment **AI**. In this case, the boundaries of **AIHG** and of **ABCD** extended, also meet at a point **J** (Figure 2).

In all cases, equation (1) holds true, i.e., a square **AIHG** has been constructed equal and similarly situated to the rectangle **ABCD**.

*Corollary 1*. *For all possible areas applied in the form of a rectangle on the fixed line segment **AB**, the point **J** traces out a parabolic curve*.

By proposition [**4**, I.45] and Problem 1 presented above, any rectilinear figure **X** can be applied to a given line segment **AB** in the simplest possible rectilinear forms of a rectangle **ABCD** and a square **AIHG** with the same area as that of **X**. Since **AB** is fixed, the areas of all possible figures **X** can be represented by a single variable, i.e., by the rectangle's height **AD** which "drags along", so to say, the side **AI** of the companion square **AIHG**. Subsequently, it seems that we have here some kind of "geometry in motion" that is best expressed by the position of the point **J** with respect to the edge-point **A** as **AD** varies. In modern terms, if the rectangle's height **AD** is the independent variable and the side **AI** of the square **AIHG** is the dependent variable, then the point **J** with respect to the edge-point **A**, geometrically represents the functional relation between these two variables. Briefly, the points **J** form a locus.

On the other hand, a similar configuration would arise if the application was effected with respect to the other edge of the line segment **AB**, the point **B**. Subsequently, the application can be considered complete in a single setting if one also took the mirror image of the preceded construction with respect to the line **AD** acting as a mirror/axis of symmetry. Let the mirror point of **J** be denoted as **J′**. Then, in the conventional language of analytic geometry, the coordinates of all points **J (x, y)** as well as those of their symmetric counterparts **J′ (−x, y)**, are related through the equation

$$x^2 = L\,y \Leftrightarrow y = \left(\frac{1}{L}\right)x^2 \qquad (2)$$

where **L** denotes the length of **AB**. This is the standard equation of the parabola in Cartesian coordinates with origin at the vertex **A** and symmetry axis the line **AD**.

Figure 3 displays three pairs of mutually symmetric **J**, **J′** points of the parabola each one of which is determined by what was demonstrated in Problem 1. The dashed lines are the tops of rectangles applied to the line segments **AB** and **B′ A** and the solid lines intersecting them are the boundaries of the companion squares.

Figure 3

Thus we have generated the parabola by means of the method of the application of areas.

# The Ellipse

**Problem 2.** *Suppose that to a given straight line segment there has been applied a rectangle equal to a given rectilinear figure and deficient by a rectangle similar to a given one. It is required to construct a square equal and similarly situated to the applied rectangle with a corner at one of the edges of the given straight line segment.*

To the given straight line **AB** there has been applied the rectangle **A B_ C_ D** equal to a given rectilinear figure **X** and deficient by a rectangle **B_ B C C_** which is similar to a given rectangular figure **Y** [**4**, VI.28]. It is required to construct a square equal and similarly situated to **A B_ C_ D** with a corner at **A**.

Let the line segment **A B_** produced by this application, be greater than or equal to **B_ C_**. Then by following the same rationale as in Problem 1, a square **AIHG** is constructed equal and similarly situated to the rectangle **A B_ C_ D** so that

$$(AB\_)(B\_ C\_) = (AG)^2 \tag{3}$$

The construction of the companion square **AIHG** is displayed in Figure 4:

Figure 4

Similarly, if the line segment **A B_** produced by this application, is less than **B_ C_**, a square **AIHG** is constructed equal and similarly situated to the rectangle **A B_ C_ D** so that equation (3) also holds true. The construction of the companion square **AIHG** is displayed in Figure 5 below:

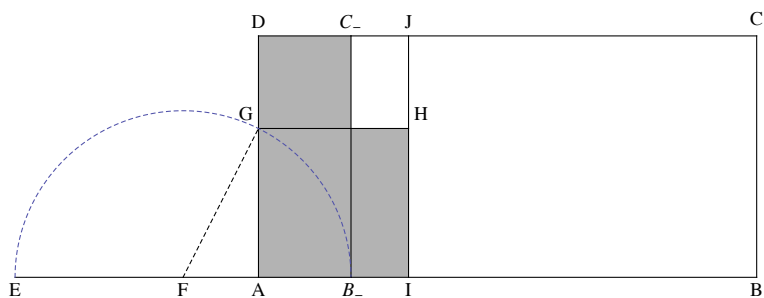

Figure 5

It is worth noting that when **A B_** is greater than or equal to **B_ C_**, the rectangle **A B_ C_ D** applies to at least half of the line segment **AB** and the boundaries of the square **AIHG** and of the rectangle **A B_ C_ D** intersect at a point **J** (Figure 4). When **A B_** is less than **B_ C_**, the rectangle **A B_ C_ D** applies to a part less than half of **AB** and the boundaries of **AIHG** and of **A B_ C_ D** extended, also meet at a point **J** (Figure 5).

In proposition [**4**, VI.27], it is demonstrated that of all rectangles **A B_ C_ D**, the one with the maximum possible area is that applied to the half of the straight line segment **AB**. Let **L** denote the length of **AB** and $\lambda$ the ratio of the sides of the rectangular figure **Y.** Then, [**4**, VI.27] reads

$$\text{AreaX} = (A\, B_-)\,(B_-\, C_-) \leq \left(\frac{L}{2}\right)\left(\frac{L}{2\lambda}\right) = \frac{L^2}{4\lambda} \tag{4}$$

This is, to our knowledge, the first ever recorded optimization problem in the history of Mathematics. As such, it offers an invaluable documented proof that *variability* constituted the hard core of the method of application of areas.

*Corollary 2*. *For all possible areas applied to a part of the fixed line segment* **AB** *in the form of a rectangle* **A B_ C_ D**, *the point* **J** *traces out an elliptic curve*.

By propositions [**4**, VI.27] and [**4**, VI.28], any rectilinear figure **X**, constrained by condition (4), can be applied to a given line segment **AB** in the form of a rectangle **A B_ C_ D** that falls short by a rectangle **B_ B C C_** similar to a given one **Y**. In Problem 2, a square **AIHG** is constructed with the same area as that of **X**. Since **AB** is fixed, the areas of all possible figures **X** can be represented by a single variable: the rectangle's height **AD**, the side **AI** of the companion square **AIHG** is also a variable that depends on **AD**, and the position of the point **J** with respect to the edge-point **A** represents geometrically the functional relation of the two. Briefly, the points **J** form a locus.

On the other hand, exactly as in Corollary 1, a similar configuration would arise if the application was effected with respect to the other edge of the line segment **AB**, point **B**. Subsequently, the application can be completed in a single setting by taking the mirror image of the preceded constructions with respect to the line **AD**. Then, the coordinates of all points **J (x, y)** and **J′ (−x, y)**, are related by

$$x^2 = L\,y - \lambda\,y^2 \tag{5}$$

This is the equation of an ellipse in Cartesian coordinates with origin at the vertex **A** and line **AD** one of the symmetry axes.

Equation (5) can be put in the standard form

$$\frac{x^2}{\left(\frac{L}{2\sqrt{\lambda}}\right)^2} + \frac{\left(y - \frac{L}{2\lambda}\right)^2}{\left(\frac{L}{2\lambda}\right)^2} = 1 \tag{6}$$

so that the basic features of the ellipse generated become manifest:

(i) The center of the ellipse is point $\left(0, \frac{L}{2\lambda}\right)$.

(ii) If $\lambda \leq 1$, the semi-major axis lies along **AD** and has a length of $\frac{L}{2\lambda}$, the semi-minor axis lies along the line parallel to **AB** that passes through the center of the ellipse and has a length of $\frac{L}{2\sqrt{\lambda}}$, and the eccentricity of the ellipse is $\sqrt{1-\lambda}$. In particular, for $\lambda = 1$, the ellipse degenerates into a circle of radius $\frac{L}{2}$.

(iii) If $\lambda > 1$, the semi-major axis lies along the line parallel to **AB** that passes through the center of the ellipse and has a length of $\frac{L}{2\sqrt{\lambda}}$, the semi-minor axis lies along **AD** and has a length of $\frac{L}{2\lambda}$, and the eccentricity of the ellipse is $\sqrt{1-\frac{1}{\lambda}}$.

Figure 6 displays three pairs of mutually symmetric **J**, **J′** points of the ellipse each one of which is determined by what was demonstrated in Problem 2. The dashed lines are the tops of rectangles applied to the line segments **AB** and **B′A** and the solid lines intersecting them are the boundaries of the companion squares.

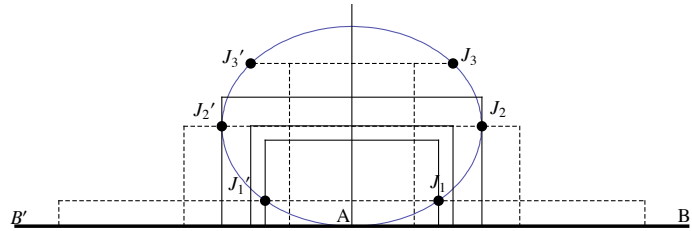

Figure 6

Thus we have generated the ellipse by means of the method of the application of areas.

## The Hyperbola

**Problem 3.** *Suppose that to a given straight line segment there has been applied a rectangle equal to a given rectilinear figure and exceeding by a rectangular figure similar to a given one. It is required to construct a square equal and similarly situated to the applied rectangle with a corner at one of the edges of the given straight line segment.*

To the given straight line **AB** there has been applied the rectangle **A B⁺ C⁺ D** equal to a given rectilinear figure **X** and exceeding by the rectangle **B B⁺ C⁺ C** which is similar to a given rectangular figure **Y** [**4**, VI.29]. It is required to construct a square equal and similarly situated to the rectangle **A B⁺ C⁺ D** with a corner at **A**.

Let the line segment **A B⁺** produced by this application, be greater than or equal to **B⁺ C⁺**. Then by following the same rationale as in Problem 1, a square **AIHG** is constructed equal and similarly situated to the rectangle **A B⁺ C⁺ D** so that

$$(AB^+)(B^+ C^+) = (AG)^2 \qquad (7)$$

The construction of the companion square **AIHG** is displayed in Figure 7:

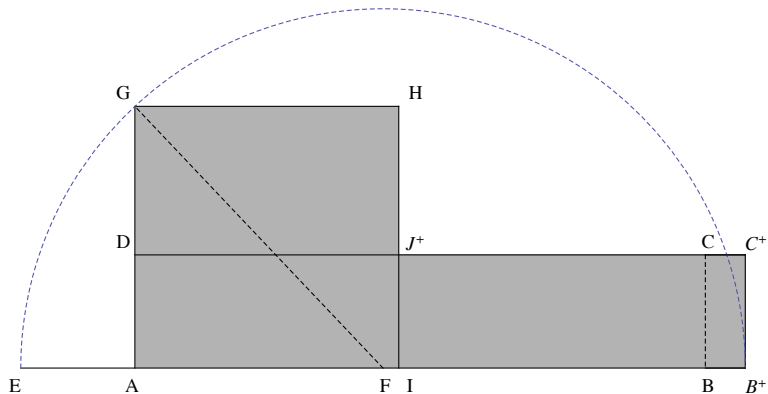

Figure 7

Similarly, if the line segment **A B⁺** produced by this application, is less than **B⁺ C⁺**, a square **AIHG** is constructed equal and similarly situated to the rectangle **A B⁺ C⁺ D** so that equation (7) also holds true. The construction of the companion square **AIHG** is displayed in Figure 8:

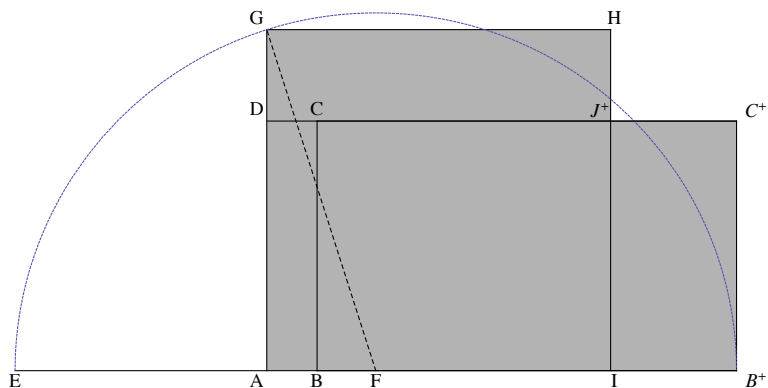

    E                        A   B   F              I            $B^+$

Figure 8

In either case, the boundaries of the constructed square **AIHG** and of the rectangle **A $B^+$ $C^+$ D** intersect at point **$J^+$**.

    ***Corollary 3***. *For all possible areas applied to the fixed line segment **AB** in the form of a rectangle **A $B^+$ $C^+$ D**, the point **$J^+$** traces out a hyperbolic curve*.

By proposition [**4**, VI.29], any rectilinear figure **X** can be applied to a given line segment **AB** in the form of a rectangle **A $B^+$ $C^+$ D** exceeding **AB** by the rectangle **B $B^+$ $C^+$ C** which is similar to a given rectangular figure **Y**. In Problem 3, a square **AIHG** is constructed with the same area as that of **X**. Since **AB** is fixed, the areas of all possible figures **X** can be represented by a single variable, i.e., by the rectangle's height **AD**, the side **AI** of the companion square **AIHG** is also a variable dependent on **AD**, and the position of the point **$J^+$** with respect to the edge-point **A**, represents geometrically the functional relation of the two. Briefly, the points **$J^+$** form a locus.

On the other hand, exactly as in Corollaries 1 and 2, a similar configuration would arise if the application was effected with respect to the other edge of the line segment **AB**, point **B**. Subsequently, the application can be completed in a single setting by taking the mirror image of the preceded constructions with respect to the line **AD**.

Let **L** denote the length of **AB** and $\lambda$ the ratio of the sides of the rectangular figure **Y**. Then, the coordinates of all points **$J^+$ (x, y)** and **$J^-$ (−x, y)**, are related by

$$x^2 = L\,y + \lambda\,y^2 \tag{8}$$

This is the equation of one branch of a hyperbola in Cartesian coordinates with origin at the vertex **A** and principal axis the line **AD**.

Equation (8) can be put in the standard form

$$-\frac{x^2}{\left(\frac{L}{2\sqrt{\lambda}}\right)^2} + \frac{\left(y + \frac{L}{2\lambda}\right)^2}{\left(\frac{L}{2\lambda}\right)^2} = 1 \tag{9}$$

so that the basic features of the hyperbola generated become manifest:

(i) The center of the hyperbola is the point $\left(0, -\frac{L}{2\lambda}\right)$.

(ii) The lines **x = 0** (i.e., **AD**) and $y = -\frac{L}{2\lambda}$ are the principal and conjugate axes of the hyperbola, respectively.

(iii) The two vertices of the hyperbola are the points **A (0, 0)** and $\mathbf{A_*}\left(0, -\frac{L}{\lambda}\right)$.

(iv) The asymptotes of the hyperbola are the lines

$$y = -\frac{L}{2\lambda} \pm \frac{1}{\sqrt{\lambda}}\,x \tag{10}$$

(v) The eccentricity of the hyperbola is $\sqrt{1+\lambda}$ .

It must be noted that equation (9) entails a second construction for Problem 3 lying symmetrically below the conjugate axis. Thus, just as the upper branch of the hyperbola is made of all points **$J^+$ (x, y)** and **$J^-$ (−x, y)**, so its lower branch is made of all points **$J_+$ (x, −y)** and **$J_-$ (−x, −y)** that also satisfy equations (8) and (9).

Figure 9 displays these four **J**-points of the hyperbola each one of which is determined by what was demonstrated in Problem 3. The dashed horizontal lines are the tops of rectangles applied in excess of the four equal line segments **AB**, **AB′**, $\mathbf{A_* B_*}$, $\mathbf{A_* B_*}'$, and the solid lines intersecting them are the boundaries of the companion squares. The dashed horizontal line passing through the midpoint of the line segment $\mathbf{AA_*}$, is the conjugate axis of the hyperbola and the inclined dashed lines are the asymptotes of the hyperbola.

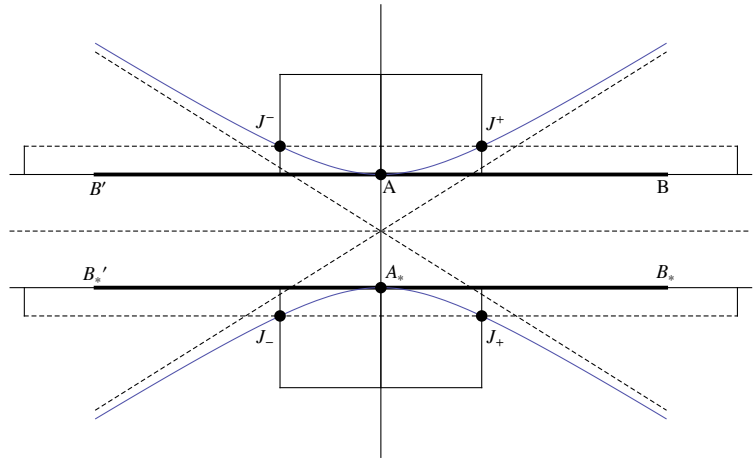

Figure 9

Thus we have generated the hyperbola by means of the method of the application of areas.

# Acknowledgments

We are indebted to our colleagues Profs A. Calogeracos and A. Boukas for their insightfull comments and suggestions.